\documentclass[12pt]{article}
\def\cleverefoptions{capitalize,noabbrev}
\usepackage[amsmath,cleveref]{e-jc}

\usepackage{graphicx}
\usepackage{makecell}
\usepackage{microtype}

\newcommand{\F}{\mathcal{F}}
\newcommand{\N}{\mathbb{N}}
\newcommand{\R}{\mathbb{R}}
\newcommand{\Q}{\mathbb{Q}}
\newcommand{\W}{\mathcal{W}}
\newcommand{\Com}{\mathbb{C}}
\newcommand{\Z}{\mathbb{Z}}
\newcommand{\B}{\mathcal{B}}
\newcommand{\s}{\mathcal{S}} % because \S is already defined
\def\one{\mathtt{1}}
\def\zero{\mathtt{0}}
\def\SEQ{\operatorname{\mathtt{SEQ}}}
\let\numeric\mathtt
\newcommand\oeis[1]{\href{https://oeis.org/#1}{#1}}

\dateline{Jan 11, 2024}{Jun 21, 2025}{Aug 22, 2025}
\Copyright{The authors. Released under the CC BY-ND license (International 4.0).}
\title{Structure and growth of $\mathbb{R}$-bonacci words}

\author{Sergey Dovgal\authornote{1}
\and
Sergey Kirgizov\authornote{2}
}

\authortext{1}{IMB UMR 5584, CNRS, Universit\'{e} Bourgogne-Europe, Dijon, France
   (\email{dovgal.alea@gmail.com}).}
\authortext{2}{LIB UR 7534, Universit\'{e} Bourgogne-Europe, Dijon, France (\email{sergey.kirgizov@u-bourgogne.fr}).}

\begin{document}

\maketitle

\begin{abstract}
  A binary word is called $q$-decreasing, for $q>0$, if inside this word
  each of length-maximal (in the local sense) occurrences of a factor of
  the form $0^a1^b$, $a>0$, satisfies $q \cdot a > b$.  We bijectively
  link $q$-decreasing words with certain prefixes of the cutting
  sequence of the line $y=qx$.  We show that for any real positive $q$
  the number of $q$-decreasing words of length $n$ grows as $C_q \cdot
  \Phi(q)^n$ for some constant $C_q$ which depends on $q$ but not on
  $n$.  From previous works, it is already known that $\Phi(1)$ is the
  golden ratio, $\Phi(2)$ is equal to the tribonacci constant, $\Phi(k)$
  is $(k+1)$-bonacci constant.  We prove that the function $\Phi(q)$ is
  strictly increasing, discontinuous at every positive rational point,
  and exhibits a fractal structure related to the Stern--Brocot tree and
  Minkowski's question mark function.
\end{abstract}

\section{Introduction}
\label{section:introduction}
\begin{figure}[ht]
  \centering
  \includegraphics[width=0.7\textwidth]{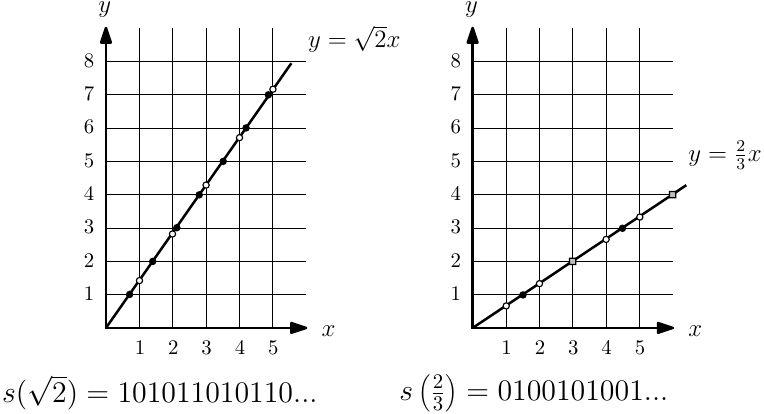}
  \caption{Cutting sequences with slopes $\sqrt{2}$ and $\frac{2}{3}$.}
  \label{f1}
\end{figure}

For any real $q>0$, the {\em ray cutting word} $s(q)$ is
defined as an intersection sequence of a straight half-line $y = qx$
for $x \in (0,\infty)$ with the lines of a square grid ($y=i$ or $x=i$ for $i
\in \N^+$). Going along the half-line, starting from $(0,0)$, we
write $\one$ if the line intersects a horizontal edge and $\zero$ in
case of a vertical edge (see~\cref{f1}), we write \( \numeric{01} \)
(in this order) when crossing an intersection point of grid lines.

For any irrational slope $q$, the word $s(q)$ is aperiodic and
Sturmian. In the general setting, Sturmian words are defined as
cutting sequences of the line $y=ax+b$ for $x\in(0,\infty)$,
irrational $a > 0$ and real $b \in [0,1)$ or equivalently as binary
  words having exactly $n+1$ factors (contiguous subwords) of length
  $n$. Sturmian words shine in several areas of mathematics:
  combinatorics, number theory, tilings, discrete dynamical systems.
  The structures similar to Sturmian words were already studied by
  Johann III Bernoulli~\cite{be} in 1771.  Expositions of Sturmian
  words and related results can be found in Chapter 2 (written by
  Berstel and S\'{e}\'{e}bold) of Lothaire's œuvre~\cite{lo} and in
  the book of Allouche and Shallit~\cite{auto}.  For a rational slope
  $q$, the word $s(q)$ is periodic, its shortest factor $f$ such that
  $s(q)=f \cdot f \cdot f \ldots$, where $\cdot$ means concatenation,
  corresponds to the Christoffel word of slope $q$~\cite{cri, ber}.

One paradigmatic example of Sturmian words is the \emph{Fibonacci
word} \( \numeric{0100101001001...} \) which is characterized by a
cutting sequence of the line with a slope \( 1/\varphi \), where \(
\varphi = (1+\sqrt 5)/2 \).  It can also be obtained either by
a recursive simultaneous application of substitution rules \( \{
\numeric{0} \mapsto \numeric{01}, \numeric{1} \mapsto \numeric{0} \}
\) to an initial word \( \numeric{0} \), or as a limit of recursive
concatenations of words \( S_n = S_{n-1}S_{n-2} \), where \( S_0 =
\zero \) and \( S_1 = \zero \one \).

Now consider another Fibonacci object, (or, more
generally, \( k \)-bonacci), which is an ensemble of binary words of
length \( n \) avoiding \( k \) consecutive \( \one \)s.
It seems that this object appears for the first time
in Knuth's book~\cite[p. 286]{knuth3}.
The set of such words is in bijection with tilings of
stripes of length \( (n+1) \times 1 \) with tiles of size \( 1 \times 1
\) (monomers), \( 2 \times 1 \) (dimers), \ldots, \( k \times 1 \)
($k$-mers), so it is convenient to call them \emph{$k$-bonacci tilings}.
The cardinality of the set of such words of length \( n \) is equal to
\( n \)th \( k \)-bonacci number (see Feinberg~\cite{fei} and Miles~\cite{miles}),
which is obtained by a recurrence relation
\( a_n = a_{n-1} + a_{n-2} + \ldots + a_{n-k} \) with initial conditions
\( a_0 = a_{-1} = 1 \) and \( a_{j} = 0 \) for any \( j < -1 \).

These two Fibonacci objects belong to two seemingly different worlds.
In this paper we propose a link between these worlds: we show that
certain subsets of prefixes of ray cutting words can be used as
building blocks to construct generalized Fibonacci tilings.  To
demonstrate this link, we extend the family of $q$-decreasing words,
defined in~\cite{ourfibo}, to cover all positive real numbers as
possible values of the parameter $q$.
\begin{definition}
  For $q \in \R^+$, a binary word is called $q$-decreasing, if inside
  this word each of length-maximal occurrences of a factor of the form
  $ \zero^a \one^b$, $a>0$, satisfies $q\cdot a > b$.  The
  length-maximality of the occurrences should be taken in the local
  sense: they are not preceded by a 0 or followed by a 1.
\label{def:1}
\end{definition}
We denote by $\W_{q,n}$ the set of {\em $q$-decreasing words} of
length $n$,~\cref{t1} gives some examples.  It is interesting to note
that E\u{g}ecio\u{g}lu and Ir\v{s}i\v{c}~\cite{ei} independently
discovered and studied hypercube subgraphs associated with a subset of
words from $\W_{1,n}$. By giving a Gray code for $\W_{1,n}$, Baril,
Kirgizov and Vajnovszki~\cite{ourfibo} prove the
E\u{g}ecio\u{g}lu-Ir\v{s}i\v{c} conjecture~\cite{ei} about the
existence of a Hamiltonian path in such hypercube subgraphs.
Recently, Wong, Liu, Lam and Im~\cite{wang} found a 2-Gray code (where
consecutive words differs in at most 2 positions) for $\W_{q,n}$ for
any real positive $q$.  The question whether there exists a 1-Gray
code when $q$ is a natural number greater than 1 remains open.

\begin{table}
  \small
  \centering
  \begin{tabular}{c|c|c|c|c|c}
    $n$ & 1 & 2 & 3 & 4 & 5\\ \hline
    $\W_{\sqrt{2},n}$
    & \makecell{
      0 \\
      1
    } &
    \makecell{
      00 \\
      01 \\
      10 \\
      11
    }
    &
    \makecell{
      000 \\
      001 \\
      010 \\
      100 \\
      101 \\
      110 \\
      111
    }
    &
    \makecell{
      0000 \\
      0001 \\
      0010 \\
      0011 \\
      0100 \\
      0101 \\
      1000 \\
      1001 \\
      1010 \\
      1100 \\
      1101 \\
      1110 \\
      1111
    }
    &
    \makecell{
      00000        \\
      00001  10010 \\
      00010  10011 \\
      00011  10100 \\
      00100  10101 \\
      00101  11000 \\
      00110  11001 \\
      01000  11010 \\
      01001  11100 \\
      01010  11101 \\
      10000  11110 \\
      10001  11111
    } \\ \hline
    $| \W_{\sqrt{2},n} | $
    & 2 & 4 & 7 & 13 & 23
  \end{tabular}
  \quad
\begin{tabular}{c|c|c|c|c|c}
  $n$ & 1 & 2 & 3 & 4 & 5\\ \hline
  $\W_{2/3,n}$
  & \makecell{
    0 \\
    1
  } &
  \makecell{
    00 \\
    10 \\
    11
  }
  &
  \makecell{
    000 \\
    001 \\
    100 \\
    110 \\
    111
  }
  &
  \makecell{
    0000 \\
    0001 \\
    0010 \\
    1000 \\
    1001 \\
    1100 \\
    1110 \\
    1111
  }
  &
  \makecell{
    00000 \\
    00001 \\
    00010 \\
    00100 \\
    10000 \\
    10001 \\
    10010 \\
    11000 \\
    11001 \\
    11100 \\
    11110 \\
    11111 \\
    \\
  } \\ \hline
  $| \W_{2/3,n} | $
  & 2 & 3 & 5 & 8 & 12
\end{tabular}
\caption{$q$-decreasing words for $n \in [1,5]$ and $q \in \{ \sqrt{2}, \frac{2}{3}
  \}$}
\label{t1}
\end{table}

In the paper~\cite{ourfibo} it has been shown that $q$-decreasing
words, for $q\in\N^+$, are in bijection with \( k \)-bonacci tilings,
where \( k=q+1 \), i.e.\  with the set of $n$-length binary words that
avoid $q+1$ consecutive \( \one \)s.  Baril, Kirgizov and Vajnovszki
are also proved that $\Phi(1)$ is the golden ratio, $\Phi(2)$ is equal
to the tribonacci constant, $\Phi(k)$ is $(k+1)$-bonacci constant.

Intriguing bijective and enumerative connections between certain
subsets of Dyck paths, integer compositions and $q$-decreasing words
are studied by Barcucci, Bernini, Bilotta and
Pinzani~\cite{barcucci_D, barcucci_R}.  Hassler, Vajnovszki and
Wong~\cite{hassler} provided, among other things, captivating Gray
codes for the set of length $n$ binary words of weight $k$ (i.e., with
exactly $k$ \( \one \)s) with the property that any prefix contains at
least $p$ times as many \( \zero \)s as \( \one \)s, for any $p \in
\R^+$.

\medskip

\begin{figure}[ht]
  \centering
\includegraphics[width=.35\textwidth]{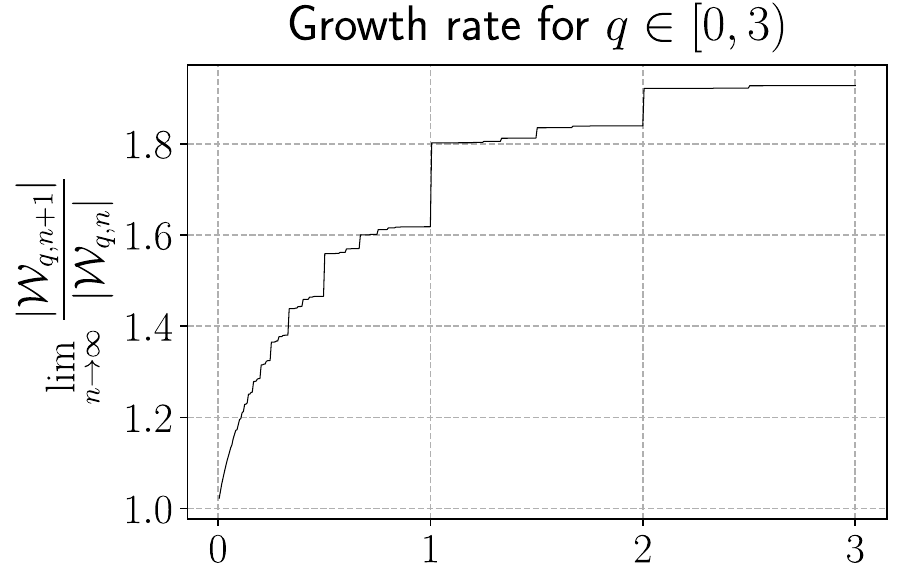}%
\includegraphics[width=.32\textwidth]{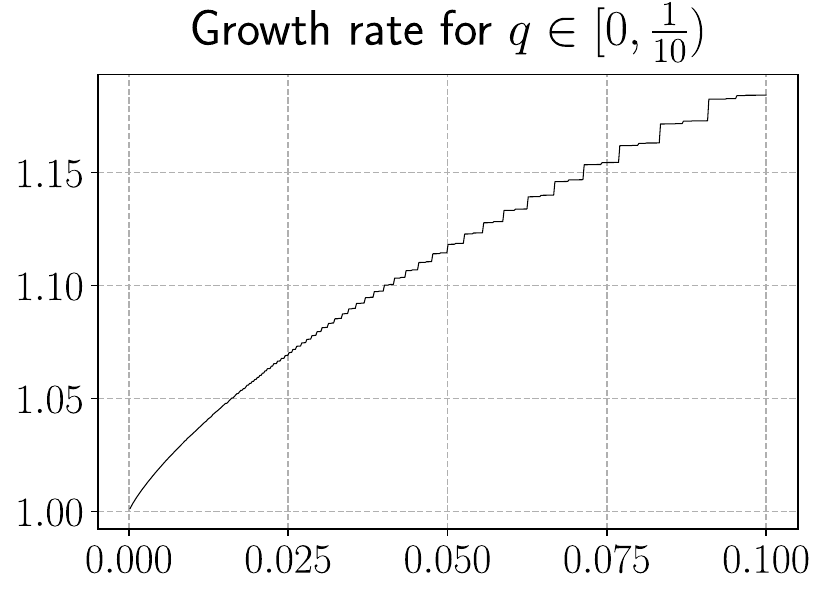}%
\includegraphics[width=.32\textwidth]{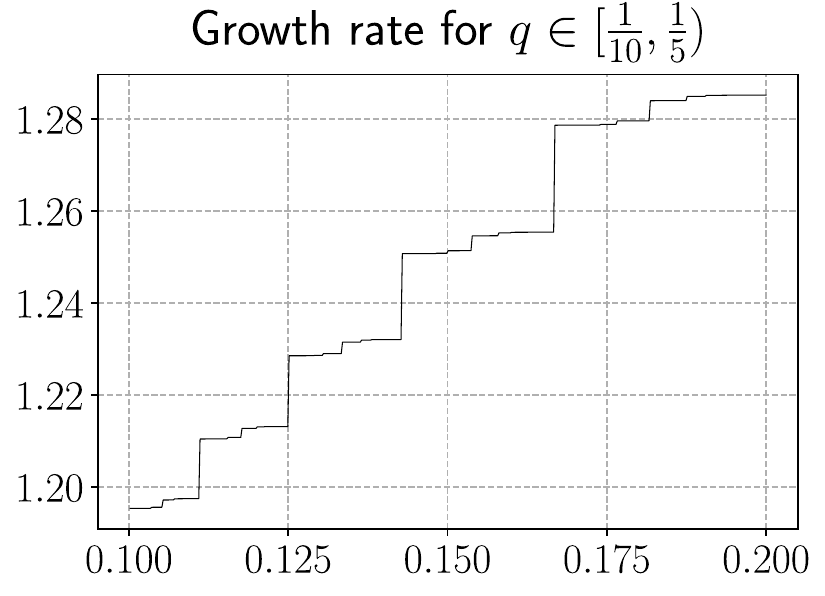}
  \caption{$\lim_{n\to\infty} |\W_{q,n+1}|/|\W_{q,n}|$ as a function of $q$
    in three different intervals.
    This function is jump discontinuous at every positive rational point.}
  \label{f2}
\end{figure}

In Section~\ref{sec2} we decompose $q$-decreasing words into sequences
of words corresponding to ray cutting prefixes ending on $1$.  The
number \( \Phi(q) \), called the exponential growth constant, is
defined as the limit ratio of successive cardinalities $\Phi(q) =
\lim_{n\to\infty} |\W_{q,n+1}|/|\W_{q,n}|$.  In
Sections~\ref{sec:radis} and~\ref{sec:thefractal} we show that this
limit exists, and explore the structure of \( \Phi(q) \) as a function
of $q$.  It turns out that the function \( \Phi(q) \) is bounded,
discontinuous at every positive rational point, strictly increasing
over \( (0, \infty) \), and also exhibits a nice fractal structure,
shown in~\cref{f2}, which bears visual resemblance to fractals arising
from information-theoretical applications such as the (appropriately
rescaled) number of coin tossings required to obtain a discrete
uniform distribution on $[1,n]$ as $n$ goes to infinity (see a
work~\cite[Fig. 6, right]{cointossing} by Bacher, Bodini, Hwang and
Tsai). A characteristic trait of such fractals is that they
demonstrate a sort of self-similarity, which is still quite tricky to
explain, as the aforementioned similarity is only approximate.  We
also show that at the vicinity of each rational point \( q \), \(
\Phi(q) \) converges locally to a piecewise linear function.

\section{Construction from ray cutting prefixes}
\label{sec2}

Here we express $q$-decreasing words as sequences of ray cutting
prefixes ending on 1. It is handy to use the Kleene star operator (it
corresponds to $\SEQ$ operator in the Flajolet--Sedgewick
book~\cite{fla}), which constructs a disjoint union of finite
concatenations from words of a given family. For instance, \( (\{
\zero, \one \zero \})^* \) provides all binary words which are empty
or end on \( \zero \) and do not contain two consecutive \( \one
\)s. We also use the ``$\cdot$'' symbol to denote all possible
pairwise concatenations between the elements of two families.

\begin{proposition}
For \( q \in \R^+ \), the set \( \W_q \) of $q$-decreasing words can be
represented as
\[
  \W_q =
  (\{ \one \})^* \cdot
  (
  \s_q
  )^*, \mbox{ where\, }
  \s_q = \cup_{i=0}^\infty
    \{ \zero^{1 + \lfloor i/q \rfloor} \one^i \}.
    \]
\label{proposition:decomposition}
\end{proposition}
\begin{proof}

By Definition~\ref{def:1}, a $q$-decreasing word is a concatenation of
factors $\zero^a \one^b$ satisfying $a=0$ or $qa > b$.  If $a=0$, the
word starts with an arbitrary sequence of $\one$s, which is \( (
\{ \one \})^*\), otherwise the condition $qa > b$ can be rewritten as $a
\geq 1 + \lfloor b/q \rfloor$. By grouping the extra zeros at the
beginning of each factor $\zero^a \one^b$, we write it as $\zero^t
\zero^{1 + \lfloor i/q \rfloor} \one^i$ with $i=b$ and
$t=b-a$. Furthermore, since $\zero \in \s_q$, the factor $\zero^t$
belongs to the family $(\s_q)^*$, and the remaining part $\zero^{1 +
  \lfloor i/q \rfloor} \one^i$ belongs to $\s_q$ itself. This
procedure allows us to decompose the remainder into a sequence of
words from $\s_q$, which finishes the proof.
\end{proof}

For a binary word $\alpha$ containing $n$ \( \zero \)s and $m$ \( \one
\)s, we define a transformation $\kappa (\alpha) = \zero^{n+1}
\one^m$, so that the empty word $\epsilon$ is mapped to the word
$\zero$. We provide a decomposition of $q$-decreasing words into
partitions of certain ray cutting prefixes by using the above
transformation.

\begin{proposition}
  For $q \in \R^+$, the transformation $\kappa$ bijectively maps the
  set of prefixes ending with 1 of the ray cutting word $s(q)$ to the
  set $\s_q$, which is used in the construction of $q$-decreasing
  words.
\end{proposition}

\begin{proof}
  Take a ray cutting word $s(q) = s_1 s_2 s_3 s_4 \ldots$ where every
  $s_i$ is a binary digit.  The index of $i$th $\one$ in this word is
  $i + \left\lfloor \frac{i}{q} \right\rfloor$.  The prefix
  $s_1 s_2 \ldots s_{i + \lfloor i/q \rfloor}$ of $s(q)$ contains exactly $i$ $\one$s and
  $\left\lfloor \frac{i}{q} \right\rfloor$ $\zero$s.  The word $\kappa(s_1
  s_2 \ldots s_{i + \lfloor i/q \rfloor}) = \zero \zero^{\left\lfloor \frac{i}{q}
    \right\rfloor} \one^i$ is a factor from the set $\s_q$, which completes the
    proof.
\end{proof}
See \cref{t2,t3} for examples.

\begin{table}[ht]
  \small
  \centering
  \begin{tabular}{c|c|c|c}
    $q$ & Ray cutting  word & Factors from $\s_q$ &
    \makecell{Some $q$-decreasing\\
      words
      }\\ \hline
    $\sqrt{2}$ &
    101011010110\ldots &
    \makecell[l]{
       $\kappa(\epsilon) = 0$, \\
       $\kappa(1) = 01$, \\
       $\kappa(101) =  0011$, \\
       $\kappa(10101) =  000111$, ...}
  &
    \makecell*{
    111100011000111,\\
    000111101000001,\\
    100000101000001}
    \\[1ex] \hline
    $\frac{2}{3}$ &
    010010100101\dots &
    \makecell[l]{
      $\kappa(\epsilon) = 0$, \\
      $\kappa(01) = 001$, \\
      $\kappa(01001) = 000011$, \\
      $\kappa(0100101) = 00000111$, ...}
    &
    \makecell*{
    111100001100001,\\
    000011001000001,\\
    001000000001111
    }
    \\ \hline
  \end{tabular}
  \caption{An illustration of the transformation $\kappa$. Prefixes
    ending with $\one$ of the ray cutting word $s(q)$ correspond to factors
    from the set $\s_q$.}
  \label{t2}
\end{table}

\begin{table}[ht]
  \scriptsize
  \centering
\begin{tabular}{c|c|c|p{5em}}
    $q$ & Ray cutting word & Counting Sequence $|\W_{q,n}|$ & OEIS \\\hline
    $\frac{1}{2}$ &
    0010010010010010... &
    1, 2, 3, 4, 6, 9, 13, 19, 28, 41, 60, 88, 129, 189, ... &
    Narayana's cows, \oeis{A930} \\ \hline
    $1/\varphi$ &
    \makecell{
      0100101001001010...
      \\
      Fibonacci word} &
    1, 2, 3, 5, 8, 12, 19, 30, 47, 74, 116, 182, 286, 448, ... &
    NEW \\ \hline
    $\frac{2}{3}$ &
    0100101001010010...
    &
    1, 2, 3, 5, 8, 12, 19, 30, 47, 74, 116, 182, 286, 449, ... &
    Comp. into 1s, 3s and 5s,
    \oeis{A60961}\\ \hline
    $1$ &
    0101010101010101... &
    1, 2, 3, 5, 8, 13, 21, 34, 55, 89, 144, 233, 377, 610, ... &
    Fibonacci, \oeis{A45} \\ \hline
    $2$ &
    1011011011011011... &
    1, 2, 4, 7, 13, 24, 44, 81, 149, 274, 504, 927, 1705, 3136, ... &
    Tribonacci, \oeis{A73} \\ \hline
    $\frac{3}{2}$ &
    1010110101101011... &
    1, 2, 4, 7, 13, 23, 42, 76, 138, 250, 453, 821, 1488, 2697, ...  &
    NEW \\ \hline
    $\sqrt{2}$ &
    1010110101101010... &
    1, 2, 4, 7, 13, 23, 42, 76, 138, 250, 453, 821, 1488, 2697, ...  &
    NEW \\ \hline
    $\varphi$ &
    1011010110110101... &
    1, 2, 4, 7, 13, 24, 44, 81, 148, 272, 499, 916, 1681, 3085, ... &
    NEW \\ \hline
    $e$ &
    1101110111011011... &
    1, 2, 4, 8, 15, 29, 56, 108, 208, 401, 773, 1490, 2872, 5536, ... &
    NEW
    \\ \hline
    $\pi$ &
    1110111011101110... &
    1, 2, 4, 8, 16, 31, 61, 120, 236, 463, 910, 1788, 3513, 6901, ... &
    NEW \\ \hline
  \end{tabular}
  \caption{Examples of ray cutting words and corresponding counting
    sequences for the cardinalities of $q$-decreasing words.}
  \label{t3}
\end{table}

\section{Rational discontinuity}
\label{sec:radis}

In this section we study the function $\Phi (q) = \lim_{n \to \infty}
\frac{|\W_{q,n+1}|}{|\W_{q,n}|}$, whose graph is shown on \cref{f2}.  Using the previously
mentioned $\SEQ$ operator, \cref{proposition:decomposition} yields the
generating function \( W_{q}(x) = \sum^\infty_{n=0} | \W_{q,n}| x^n \)
of the family \( \s_q \) for any \( q \in \R^+ \):
\begin{equation}
\label{eq2}
  W_q(x) =
  \frac{1}{
    (1-x) \left(
      1 - \sum_{i=0}^\infty
        x^{1 + i + \left\lfloor \frac{i}{q} \right\rfloor}
    \right)
  }.
\end{equation}
The case where $q$ is a positive \emph{rational number} represented by an
irreducible
fraction $\frac{c}{d}$ is treated
in~\cite{ki} where the author expresses the generating function $W_{q}(x)$ as
\begin{equation}
W_{q=\frac{c}{d}}(x) = \frac{1 -
  x^{c+d}}{(1-x)\left(1-x^{c+d} - \sum_{i=0}^{c-1} x^{1 + i +
  \left\lfloor \frac{i}{q} \right\rfloor} \right)}.
\label{eq1}
\end{equation}
In other words, \cref{eq2} holds for any positive real $q$, and is
more general, although simpler, form of \cref{eq1} which is only
valid for $q\in \Q^+$.

To prove the results about the asymptotic behaviour of coefficients
$[x^n]W_q(x)$, we need the following lemma, which can be considered as
a simpler variant of Daffodil Lemma from Flajolet--Sedgewick
book~\cite[Lemma IV.1, p. 266]{fla}.

\begin{lemma}[``Little Narcissus Lemma'']

  Let $f(x) = x + \sum_{k=2}^\infty a_k x^k$ be a power series with
  non-negative real coefficients $\{a_k\}_{k=2}^\infty$ and a positive radius of
  convergence $\rho$.  The following hold:\\
  \hspace*{2em}{\bf 1)} there is a unique real root $R$ of the
  equation $f(x)=1$;\\
  \hspace*{2em}{\bf 2)} $R$ has the multiplicity one;\\
  \hspace*{2em}{\bf 3)} there is no other root $x \in \Com, x \neq R$ of the
  equation $f(x)=1$ such that $|x| \le R$.
  \label{little-narcissus-lemma}
\end{lemma}

\begin{figure}[ht]
  \includegraphics[width=\textwidth]{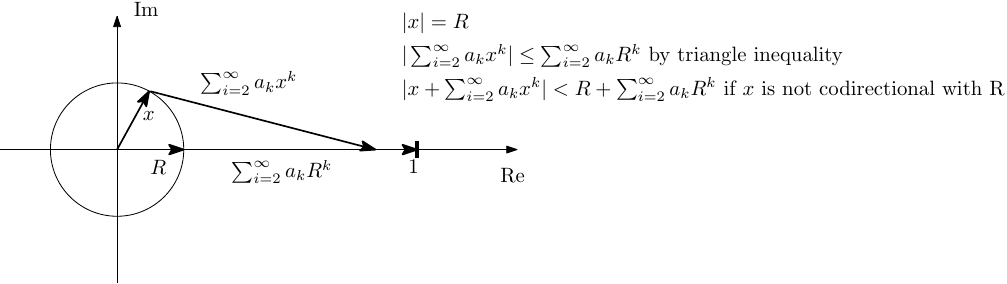}
  \caption{An element of the Little Narcissus Lemma proof.}
  \label{fig:little-narcissus}
\end{figure}

\begin{proof}
 Suppose that at least one of the coefficients $\{a_k\}_{k=2}^\infty$
 is strictly positive, otherwise the result is trivially simple.  The
 function $f(x)$ is monotonically increasing on $[0,\rho)$.  We have
   $f(0) = 0$, $\lim_{x\to \rho} f(x) = +\infty$. So, there is a
   unique real value $R \in [0,\rho)$ such that $f(R) = 1$.  The
   multiplicity of $R$ is one because $f'(R) > 0$.  For any $x$ such
   that $|x| < R$, we have, by the triangle inequality, $|f(x)| \leq
   f(|x|)<f(R)=1$.  Consider the case $|x| = R$.  Again, by the
   triangle inequality we have $|\sum_{i=2}^\infty a_k x^k| \le
   \sum_{i=2}^\infty a_k R^k$.  If $x \ne R$, we have a strict
   inequality $|x + \sum_{i=2}^\infty a_k x^k| < R + \sum_{i=2}^\infty
   a_k R^k$, see~\cref{fig:little-narcissus}.
\end{proof}

We have the following result about the asymptotic behaviour of $[x^n] W_{q}(x)$.
\begin{proposition}
  For any real $q > 0$, the number of $q$-decreasing words of length $n$ grows as $C_q \cdot
  \Phi(q)^n$, where $1/\Phi(q)$ is the unique smallest in modulus root
  of $1 - \sum_{i=0}^\infty x^{1+i+ \left\lfloor \frac{i}{q}
    \right\rfloor},$ and $$C_q = - \frac{\Phi(q)}{ \left(
    (1-x) \left(
    1 - \sum_{i=0}^\infty x^{1+i+ \left\lfloor \frac{i}{q} \right\rfloor}
    \right)
    \right)' \big( 1/\Phi(q) \big) }.$$
  \label{propCq}
\end{proposition}

\begin{proof}
  Let $f_q(x)=\sum_{i=0}^\infty x^{1 + i + \left\lfloor \frac{i}{q}
    \right\rfloor} = x + x^{2 + \left\lfloor \frac{1}{q}
    \right\rfloor} + \ldots$.  For any $x\in[0,1)$ we have
    $\sum_{i=0}^\infty x^{1 + i + \left\lfloor \frac{i}{q}
      \right\rfloor} \le \sum_{i=1}^\infty x^i \le \frac{x}{1-x}$.
    So, $f_q(x)$ evaluated at $x \in [0,1)$ is bounded,
      thus convergent for any real $q \in (0,\infty)$.
    Apply Little Narcissus Lemma to see that $f_q(x)=1$ has a unique
    smallest in modulus root $R < 1$ which has the multiplicity 1.
    Let $g_q(x) = (1-x) \left( 1 - \sum_{i=0}^\infty x^{1 + i +
      \left\lfloor \frac{i}{q} \right\rfloor} \right)$.  The function
    $W_q(x) = 1/g_q(x)$ is meromorphic in the unit disc. The root $\rho_q$ of
    $f_q(x)=1$ is the unique smallest in modulus pole of $W_q(x)$. The
    pole has the multiplicity 1.  Let $\Phi(q) = 1/\rho_q$. Using
    classical asymptotic analysis of meromorphic functions (see
    Flajolet--Sedgewick book~\cite{fla}, Sedgewick's online
    course~\cite{slides-sedgewick} or Orlov's paper~\cite{orlov}) we
    obtain the result, expressing $C_q$ as $\Phi(q) / g'(1/\Phi(q))$.
\end{proof}

\begin{figure}[ht]
  \centering
  \includegraphics[width=0.75\textwidth]{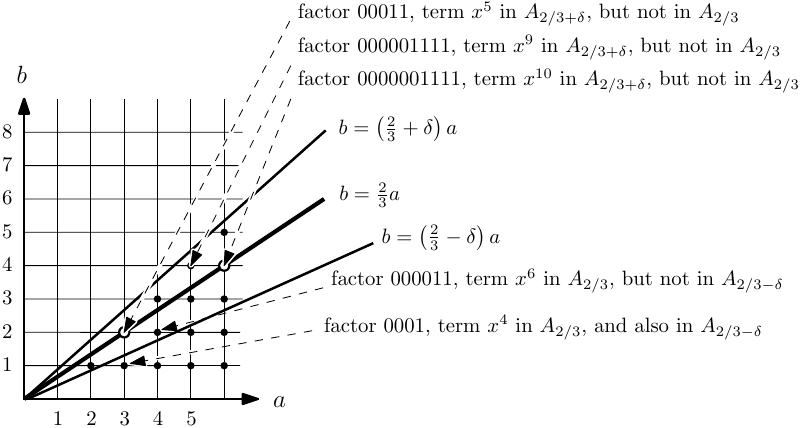}
  \caption{The line $b=\frac{2}{3}a$ and the geometrical
    interpretation of factors $0^a1^b$ where $\frac{2}{3}a>b$.}
  \label{f3}
\end{figure}

The equation $
1 = \sum_{i=0}^\infty x^{1+i+ \left\lfloor
  \frac{i}{q} \right\rfloor}
$
shares the smallest in modulus root with 
$$ A_q := 1 - \sum_{i=1}^\infty \sum_{j=0}^\infty x^{1+i+ \left\lfloor
  \frac{i}{q} \right\rfloor + j} ,
$$
a fact that has a nice geometrical interpretation.  To see it, we
have
to decompose the set $\W_q$ in another way, different from what was
given in \cref{sec2}. Here we use a set $\F_q$ of factors
$0^a1^b$ such that $qa>b$ and $b\ge1$.  With this, we decompose any
word $w\in \W_q$ as a sequence of \( \one \)s, followed by a sequence of
factors from $\F_q$, followed by a sequence of \( \zero \)s. Any of these
sequences can be empty. We have
$$
  \begin{array}{c}
    w = \overbrace{\one...\one}^{\text{some ones}}
    \overbrace{f_1 f_2 .... f_k}^{f_\ell \in \F_q}
    \overbrace{\zero...\zero}^{\text{some zeros}}, \quad
    \text{where }
    \F_q = \bigcup_{i=1}^\infty \bigcup_{j=0}^\infty \{
\overbrace{\zero\ldots\numeric{00}}^{1 + \left\lfloor \frac{i}{q} \right\rfloor + j
  \text{ zeros}}\mkern-13mu \underbrace{\one\ldots\numeric{11}}_{i \text{ ones}}
\;\; \}.
  \end{array}
$$
Now, we write the g.f. $W_q(x)$ as
$$
  W_q(x) = \frac{1}{1-x} \cdot
  \frac{1}{A_q}
  \cdot \frac{1}{1-x}.
$$
Consider the grid $\Z^+ \times \Z^+$, and make every point $(a,b)$
correspond to a factor $\zero^a \one^b$.  The power series $A_q$ sums over all
points with positive integer coordinates found under the line
$b=qa$. \cref{f3} gives some examples.

From \cref{propCq} we see that the function $\Phi (q)$ is well-defined
as $\lim_{n \to \infty} \frac{|\W_{q,n+1}|}{|\W_{q,n}|}$.  In
\cref{propPhi} we prove the basic properties of this strictly
increasing and bounded function with a countable number of
discontinuities, and then calculate the the jump sizes in
\cref{proposition:jump}.
  
\begin{proposition}
  The function
  $ \Phi (q) = \lim_{n \to \infty} \frac{|\W_{q,n+1}|}{|\W_{q,n}|}$ is\\
  \hspace*{2em}{\bf a)} strictly increasing over $q \in [0, \infty)$;\\
  \hspace*{2em}{\bf b)} bounded, $1 \le \Phi(q) < 2$, with $\Phi(0)=1$ and
  $\lim_{q\to\infty}\Phi(q)=2$; \\
  \hspace*{2em}{\bf c)} left-continuous (and right-discontinuous) at every positive rational point;\\
  \hspace*{2em}{\bf d)} continuous at every positive irrational point.
  \label{propPhi}
\end{proposition}
\begin{proof}

  Let $\delta > 0$.

  {\bf a)} Note that $ A_{q+\delta} = 1 - \sum_{i=1}^\infty
  \sum_{j=0}^\infty x^{1+i+ \left\lfloor \frac{i}{q+\delta}
    \right\rfloor + j} $ contains all terms of $ A_q = 1 -
  \sum_{i=1}^\infty \sum_{j=0}^\infty x^{1+i+ \left\lfloor \frac{i}{q}
    \right\rfloor + j} $ together with terms not presented in $A_q$.
  This can be seen geometrically on \cref{f3}. The terms corresponding
  to points of positive integer coordinates, lying at the line $b=qa$,
  are included in $A_{q+\delta}$ for any $\delta > 0$, but none of
  these terms are in $A_{q}$. So, we necessarily have $\rho_{q+\delta}
  < \rho_q$, and thus $\Phi(q) < \Phi(q+\delta)$.

  {\bf b)} If $q=0$, the only $q$-decreasing word of length $n$ is
  $\one^n$, so $\Phi(0) = 1$. It is straightforward that $\lim_{q
    \to 0} \Phi(q) = 1$.  As $q \to \infty$ we allow more and more
  binary words, and the functional limit of $A_q$ can be expressed as 1 minus
  the sum over all integer points from the positive quadrant:
  \(
    \lim_{q \to \infty} A_q(x) = 1-\sum_{i=1}^\infty \sum_{j=1}^\infty x^{i+j}
  \),
  so $\lim_{q \to \infty} \Phi(q) = 2$.

  {\bf c)} For a positive rational $q$ represented by an irreducible
  fraction $\frac{c}{d}$ the line $b=qa$ contains points of integer
  coordinates $(kd,kc)$ for $k \in [1, \infty)$.  For any $\delta>0$,
    these points are below the line $b=(q+\delta)a$.  The
    corresponding terms of the form $x^{kd+kc}$ are included in
    $A_{q+\delta}$ but not in $A_q$.  This, in turn, influences the
    smallest in modulus root of $A_{q+\delta}$ which is strictly less
    than the smallest in modulus root of $A_q$ and we obtain $\Phi(q)
    < \lim_{\delta \to 0^+} \Phi(q+\delta)$.  The difference between
    $\lim_{\delta \to 0^+}$ and $\Phi(q)$ is explicitly calculated
    in~\cref{proposition:jump}. In other case, no line of the form
    $b=(q-\delta)a$ can have the points $(kd,kc)$ below it (even if
    $\delta = 0$).  Any point lying under $b=qa$ also lies under
    $b=(q-\delta)a$ for a sufficiently small $\delta$. We obtain
    $\lim_{\delta \to 0^+} \Phi(q-\delta) = \Phi(q)$.

  {\bf d)} For a positive irrational $q$, the line $b = q a$ contains
  no points of positive integer coordinates.  There are therefore no
  terms included in $A_q$ corresponding to this line.  So, the
  smallest root $\rho_q$ of $A_q$ can be approached by smallest roots
  of $\{A_{r_i}\}_{i=1}^\infty$ where $\{r_i\}_{i=1}^\infty$ is a
  sequence of rational numbers such that $\lim_{i \to \infty } r_i =
  q$.
\end{proof}

Assume that $q\in\Q^+$ is represented by an irreducible fraction
$c/d$.  From \cref{eq1} we see that the growth rate is dictated by the
smallest in modulus root $\rho_q$, of the polynomial
$$
\Pi_q := 1 - x^{c+d} - \sum_{i=0}^{c-1} x^{1 + i + \left\lfloor
  \frac{i}{q} \right\rfloor}.
$$
Comparing \cref{eq1,eq2} we see that $\Pi_q$ shares the same smallest
in modulus root with
$$
1 - \sum_{i=0}^\infty x^{1+i+ \left\lfloor
  \frac{i}{q} \right\rfloor}.
$$

\begin{lemma}
  The smallest in modulus root $\rho_q$ of $\Pi_q$ is positive and real.
  \label{small-lemma}
\end{lemma}
\begin{proof}
  By Little Narcissus Lemma~\ref{little-narcissus-lemma}.
\end{proof}

\begin{proposition}
  For any irreducible $q = \frac{c}{d} \in \Q^+$, $ \lim_{\delta\to 0^+}
  \Phi(q+\delta) $ equals the reciprocal of the smallest in modulus root, denoted
  by $\rho^+_q$, of the polynomial 
  $$
  \Pi^+_q := 1 - (2-x) x^{c+d} -
  \sum_{i=0}^{c-1} x^{1 + i + \left\lfloor \frac{i}{q}
    \right\rfloor}.
  $$
  We have
  $$ \lim_{\delta\to 0^+} \Phi(q+\delta) - \Phi(q) =
  \frac{1}{\rho^+_q}
  - \frac{1}{\rho_q}
  ,$$ where $\rho_q$ is the smallest in modulus root of $\Pi_q
  := 1 - x^{c+d} - \sum_{i=0}^{c-1} x^{1 + i + \left\lfloor
    \frac{i}{q} \right\rfloor}.$
  \label{proposition:jump}
\end{proposition}

\begin{proof}
  Assume that $q$ is represented by an irreducible fraction $c/d$. We
  consider a set of binary words $\W^+_{q,n}$ where every
  length-maximal (in the local sense) occurrence of a factor of the
  form $\zero^a \one^b$ respects $aq \ge b$. It differs from the
  definition of $q$-decreasing words which encloses a strict
  inequality.  It is clear (see e.g. \cref{f3}) that the set
  $\W_{q+\delta,n}$ approaches $\W^+_{q, n}$ as $\delta \to 0$,
  factors $\zero^{kd}\one^{kc}$ corresponding to the points of the
  form $(kd,kc)$ are included in both sets, no points located strictly
  above the line $b=qa$ are considered in both cases.  The set
  $\W^+_q$ is constructed as $\W^+_q = (\{\one\})^* \cdot (\s_q^+)^*$,
  where $\s^+_q = \{\zero\} \cup \bigcup_{i=1}^\infty \{ \zero^{
    \lceil i/q \rceil } \one^i\}$
  (c.f. \cref{proposition:decomposition}).

  Note that $\s^+_q = \{ \zero \} \cup \B$, where $\B = \B_0 \cup \B_1 \cup
  \B_2 \cup ... $, and $\B_0 = \bigcup^{c-1}_{i=1}\{ \zero^{1+\lfloor i/q
    \rfloor} \one^i \} \cup \{ \zero^d \one^c\}$, and $\B_{j+1}$ is constructed by
  inserting the factor $\zero^d \one^c$ after the last $\zero$ in words from $\B_j$.
  So, the g.f. $W^+_q(x)$ of the words in $\W^+_q$ is 
  $$
  W^+_{q=\frac{c}{d}}(x) =
  \frac{1}{(1-x)\left( 1 - \left(x + \frac{\sum^{c-1}_{i=1} x^{1+i+\left\lfloor \frac{i}{q} \right\rfloor} + x^{c+d} }{1-x^{c+d}} \right) \right)}
  =$$
  $$
  \hspace*{5.5em} = \frac{1-x^{c+d}}{(1-x)\left( 1 + x^{c+d+1}  - 2x^{c+d} - \sum^{c-1}_{i=0} x^{1+i+\left\lfloor \frac{i}{q} \right\rfloor} \right) }
  = \frac{1-x^{c+d}}{(1-x) \cdot \Pi^+_q}.
  $$ Consider $f_q(x) = x + \frac{\sum^{c-1}_{i=1} x^{1+i+\left\lfloor
      \frac{i}{q} \right\rfloor} + x^{c+d} }{1-x^{c+d}}$, from Little
  Narcissus Lemma~\ref{little-narcissus-lemma} it follows that there
  are no complex roots of the equation $f_q(x) = 1$ smaller than or
  equal in modulus to $\rho_q^+$ that lies in $[0,1)$.  The claimed
    result is obtained by comparing this formula with \cref{eq1}.
\end{proof}

\section{The fractal}
\label{sec:thefractal}

As we can see on \cref{f2} the graph of the function $\Phi(q)$ shows a
certain amount of self-similarity. We explain some aspects of this
fractality in this section. Firstly, we zoom into the intervals $q \in
(\frac{k}{k+1}, 1]$ for $k \in [1, \infty)$, and then look into a more
    general setting using a rescaling based on Minkowski's question
    mark function and the Stern--Brocot tree.  We use this tree to
    generate a sequence of nested intervals, narrower each time,
    around a given rational number.  The summary of results and one
    open question are presented at the end of this section.
    
\subsection{Around the point $q=1$}
    
Let us recall that the smallest in modulus root of $1 - x - x^2$ is
$\rho_1 = 1/\varphi$, where $\varphi = (1+\sqrt 5)/2$.  The
polynomials \( \Pi_q \) and \( \Pi^+_q \) are defined
in~\cref{proposition:jump}.

\begin{proposition}
\label{proposition:asymptotics:1plus}
  For natural $k \ge 1$, the smallest in modulus root $\rho_{k/(k+1)}$
  of $\Pi_{k/(k+1)}$ is
  $$
  \rho_1 + C \rho_1^{2k} (1+ o(1) ),
  \quad \text{as} \quad k \to \infty,
  $$
  where $C=\rho_1^3/(1+ 2 \rho_1)$, and $\rho_1$ is the smallest in modulus root of $1-x-x^2$.
  \label{prop4}
\end{proposition}

\begin{proof}
  After arithmetic transformations, the equation
  $\Pi_{k/(k+1)} = 0$, i.e.\ the equation
  $$
  1  - x^{2k+1} - \sum^{k-1}_{i=0} x^{1+i+\left\lfloor \frac{i(k+1)}{k} \right\rfloor} = 0
  $$
  turns into
  \begin{equation}
    \nonumber
      1  - x^{2k+1} - \sum^{k-1}_{i=0} x^{1+2i} = 0
      \quad \Leftrightarrow \quad
      1 = x \cdot \frac{x^{2k+2} - 1}{x^2-1}
      \quad \Rightarrow \quad
      1 = x + x^2 - x^{2k+3}.
  \end{equation}
  We multiplied by $(x^2 - 1)$ both sides of equation, adding two new
  roots $1$ and $-1$. This does not change the overall picture of the
  asymptotics, because $0 < \rho_{k/(k+1)} < 1$.
  By~\cref{small-lemma} there is no complex roots smaller than or equal in modulus to $\rho_{k/(k+1)}$.
  The root can be represented
  as $\rho_{k/(k+1)} = \rho_1 + \varepsilon_k$ for some positive
  $\varepsilon_k$ such that $0 < \rho_1 + \varepsilon_k < 1$.  Note that
  $\varepsilon_k \to 0$ as $k$ grows.  Next, we substitute $\rho_1 +
  \varepsilon_k $ for $x$ in $x+x^2-1 = x^{2k+3}$, use $1 = \rho_1 +
  \rho_1^2$ and obtain the following: { \everymath={\displaystyle}
    $$
    \begin{array}{l}
       \rho_1 + \varepsilon_k + (\rho_1 + \varepsilon_k)^2 - 1 = (\rho_1 + \varepsilon_k)^{2k+3}, \\
       \rho_1 + \varepsilon_k + \rho_1^2  + 2 \rho_1 \varepsilon_k + \varepsilon_k^2 -1 = \rho_1^{2k+3} (1 + o(1)),  \\
      \varepsilon_k ( 1  + 2 \rho_1 + \varepsilon_k) = \rho_1^{2k+3} (1 + o(1)).  \\
    \end{array}
    $$ }%
  The claimed result $\rho_{k/(k+1)} = \rho_1 + C \rho_1^{2k}
  (1+ o(1) )$ follows, because $1 + 2 \rho_1 + \varepsilon_k \to 1 + 2
  \rho_1$ as $k \to \infty$.
\end{proof}

Using the Taylor expansion several times one can improve the root
approximation and get $\rho_{k/(k+1)} = \rho_1 + C \rho_1^{2k} +
O(k\rho_1^{4k})$. But in context of this paper,
\cref{prop4} is sufficient.

\begin{proposition}
\label{proposition:asymptotics:1minus}
  For natural $k \ge 2$, the smallest in modulus root
  $\rho^+_{(k-1)/k}$ of $\Pi^+_{(k-1)/k}$ is $$ \rho_1 + C \rho_1^{2k}
  (1+ o(1) ),$$ where $C=\rho_1^3/(1+ 2 \rho_1)$, and $\rho_1$ is the smallest in modulus root of $1-x-x^2$.
  \label{prop5}
\end{proposition}

\begin{proof}
After arithmetic transformations similar to the ones in
\cref{proposition:asymptotics:1plus}, the equation
$\Pi^+_{(k-1)/k} = 0$, i.e.\ the equation
$$
1  - x^{2k-1} - \sum^{k-2}_{i=0} x^{1+i+\left\lfloor \frac{ik}{k-1} \right\rfloor}  - x^{2k-1} + x^{2k} = 0
$$
turns into
$$x + x^2 - 1 = x^{2(k-1)} (2 x^3 - x - x^4 + x^2).$$
  Note that, at some point, we multiply both sides of equation by
  $(1-x^{2})$, adding $-1$ and $1$ as roots.
  This does not change the picture dramatically, because $0 < \rho^+_{(k-1)/k} < 1$.

  As in the proof of \cref{prop4} the root is $\rho_1 + \varepsilon_k$, $\varepsilon_k \to 0$ when $k\to\infty$.
  Note that $2 x^3 - x - x^4 + x^2 - x^5 = (x^3 - x)(1 - x - x^2)$, so we have
  $2 \rho_1^3 - \rho_1 - \rho_1^4 + \rho_1^2 = \rho_1^5$, so $2 (\rho_1 + \varepsilon_k)^3 - (\rho_1 + \varepsilon_k) -
  (\rho_1 + \varepsilon_k)^4 + (\rho_1 + \varepsilon_k)^2 = \rho_1^5 (1 + o(1)) $.
  Using the techniques from the previous proof we obtain the claimed result.
\end{proof}

\begin{figure}[ht]
  \centering
  \includegraphics[width=\textwidth]{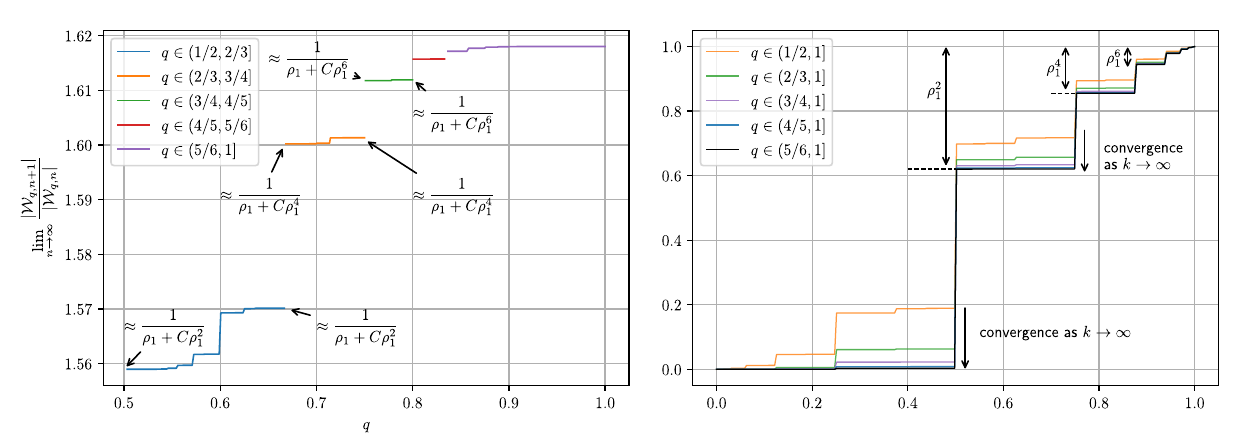}
  \caption{Fractal structure of the function
    $ \Phi(q) = \lim_{n\to\infty} \frac{|\W_{q,n+1}|}{|\W_{q,n}|}$,
    before and after rescaling on the intervals $(\frac{k-1}{k}, 1]$.
  }
  \label{f4}
\end{figure}

\subsection{Minkowski rescaling and Stern--Brocot tree}

Note that
\cref{proposition:asymptotics:1plus,proposition:asymptotics:1minus}
already provide some insights into the fractal structure of $\Phi(q)$
displayed in \cref{f4}. On the left side of this figure, the
horizontal axis is partitioned into intervals $(1/2, 2/3], (2/3, 3/4],
  \ldots$, and the parts of the plot of $\Phi(\cdot)$ are grouped
  accordingly. In particular, we showed that
$$
    \lim_{k \to \infty} \dfrac{\Phi(1) - \lim_{\delta \to 0^+}\Phi(\frac{k-1}{k} + \delta)}
    {\Phi(1) - \Phi(\frac{k}{k+1})} = 1,
$$
i.e.\ that the images of the intervals
$(\frac{k-1}{k}, \frac{k}{k+1}]$ tend to straighten as $k \to \infty$.

To better observe the self-similarity, the intervals $q \in (k/(k+1), 1]$
and their images under $\Phi(\cdot)$
can be ``normalized'' using {\em simple rescaling} and {\em
    Minkowski's question-mark function}~\cite{denjoy, minkowski}.
  Simple rescaling takes a set of positive values $V$, containing at
  least 2 values, and maps every $v \in V$ to $\frac{v - \min V}{\max V
    - \min V}$, so the image lies in $[0,1]$.  Minkowski's
  question-mark function is a little trickier, and we must first
  discuss mediants and the construction of the Stern--Brocot
  tree~\cite{stern, brocot}.

\begin{figure}[ht]
  \centering
  \includegraphics[width=0.7\textwidth]{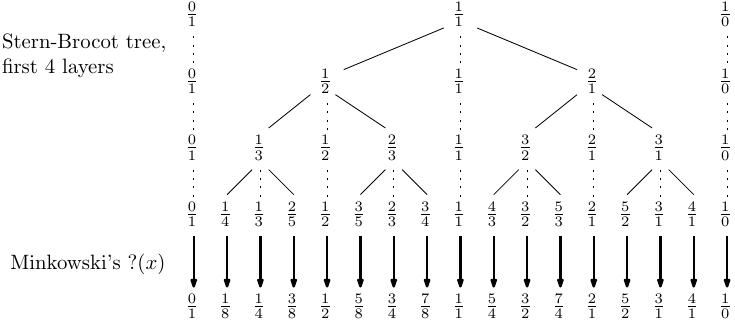}
  \caption{The Stern--Brocot tree and Minkowski's $?(x)$.}
  \label{f5}
\end{figure}

For two irreducible fractions $a/b$ and $c/d$ their {\em mediant} is
defined as $(a+c)/(b+d)$. The root of the Stern--Brocot tree is $1/1$,
which is the mediant of two conventionally irreducible fractions $1/0$
and $0/1$. To determine the left (resp. right) child of a node $x/y$
of the level $i$ we need to find the greatest (resp. smallest)
fraction $x'/y' < x/y$ (resp. $x'/y' > x/y$) that appears in set of
values of first $i$ levels together with $1/0$ and $0/1$, and compute
the mediant $(x+x')/(y+y')$. For instance, the left child of $2/3$ is
$3/5$, it is calculated as the mediant of $1/2$ and $2/3$. \cref{f5}
illustrates this process.  Stern--Brocot tree contains all rationals
once.

Minkowski's question-mark function, denoted by $?(x)$, maps a positive
rational value $x$ to a positive dyadic rational $a/2^k$ with $a,k\in
\N$. By definition, $?(0) = 0$ and $?(1) = 1$.  Whenever $x \in (0,1)$
is a rational number represented by an irreducible fraction $a/b$,
such that in the Stern--Brocot tree it is constructed via taking a
mediant of two fractions $p/r$ and $p'/r'$, its image under
Minkowski's function is defined as
$$
    \mathop{?}\left(x\right)
    = \mathop{?}\left(
        {\frac {p+p'}{r+r'}}
    \right)
    :={\frac {1}{2}}\left(
        \mathop{?} \left({\frac {p}{r}}\right)
        +\mathop{?} \left({\frac {p'}{r'}}\right)
    \right).
$$
In other words, we descend the Stern--Brocot tree in search of the
$a/b$, and ``in parallel'' construct a resulting value by applying the
mean instead of the mediant. For $x>1$, Minkowski's function is
defined as $?(x+1) = ?(x) + 1$. In general, $?(x)$ is monotonically increasing,
and can be defined on all $\R^+$~\cite{denjoy}.

\begin{figure}[ht]
  \centering
  \includegraphics[width=.55\textwidth]{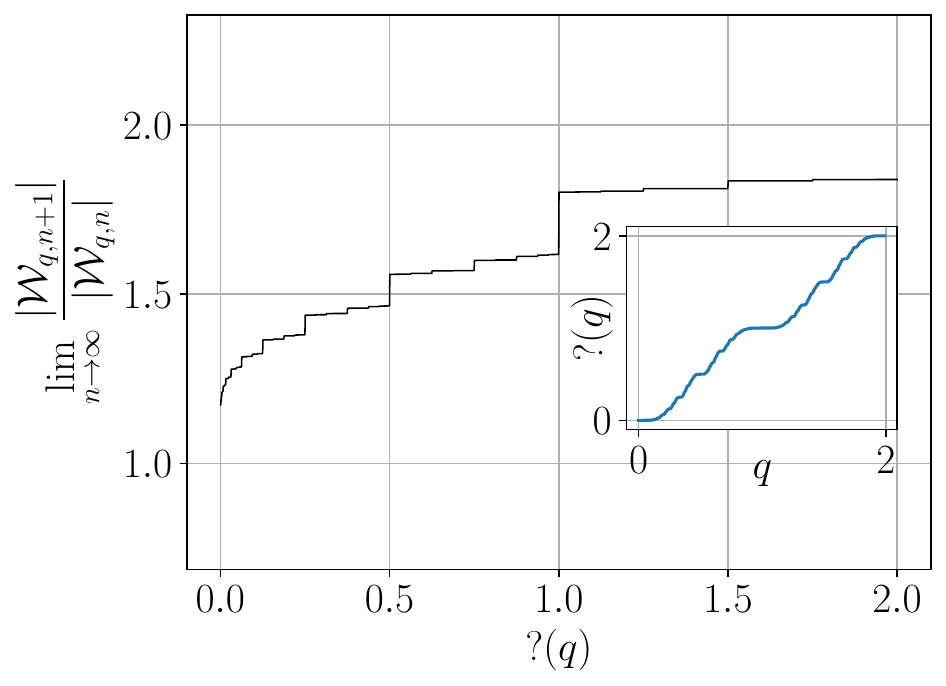}
  \caption{$\Phi(q) := \lim_{n\to\infty} |\W_{q,n+1}|/|\W_{q,n}|$ as a function of $?(q)$.}
  \label{f7}
\end{figure}

The right side of \cref{f4} is obtained by applying the simple
rescaling on the vertical axis and Minkowski's question-mark function
followed by the simple rescaling on the horizontal axis for intervals
$(k/(k+1), 1]$ and their images. The similar analysis can be done for
   intervals $(1, (k+1)/k]$. The fractal structure of $\Phi$ presented
 in~\cref{f2} appears more regular in~\cref{f7} as we apply
 Minkowski's question-mark function over the $x$-axis.

\clearpage
\subsection{Around a positive rational number $q$}

Now, we study the more general case, that is the fractal structure of
$\Phi(q)$ on the intervals $\left(\frac{p+ck}{r+dk},\frac{c}{d}
\right]$ for $k \ge 1$, where fractions are irreducible, and
$\frac{p+c}{r+d}$ is the left child of $c/d$ in the Stern--Brocot
tree.  Next proposition can be applied, for example, to the intervals
$\left( \frac{3+5k}{2+3k},\frac{5}{3} \right]$, where $8/5$ is the
  left child of $5/3$.

\begin{proposition}
  For natural $k \ge 1$, the smallest in modulus root of $\Pi_{\frac{p+ck}{r+dk}}$ is
  $$\rho_{\frac{p+ck}{r+dk}} = \rho_{c/d} + C 
  \rho^{(c+d)k}_{c/d} \big( 1 + o(1) \big),$$ where $k\ge 1$,
  $\frac{p+c}{r+d}$ is the left child of $c/d$ in the Stern--Brocot tree,
  $C$ is a constant depending only on $p/r$ and $c/d$, and $\rho_{c/d}$ is the
  smallest in modulus root of $\Pi_{c/d}$.
  \label{prop6}
\end{proposition}

\begin{proof}
   Recall that from~\cref{proposition:jump} we have
  \begin{equation}
    \Pi_{\frac{p+ck}{r+dk}}:=
  1 - x^{p+r+(c+d)k} - \sum_{i=0}^{p+ck-1} x^{1 + i + \left\lfloor
    \frac{i (r+dk) }{p+ck} \right\rfloor}.
  \label{eq3}
  \end{equation}
  For $0<i<p+ck$, $1 + i + \left\lfloor \frac{i (r+dk) }{p+ck}
  \right\rfloor$ equals the number of integer points with coordinates
  $(i, y)$ lying between two diagonal lines intersecting at the origin
  with respective slopes $(-c)$ and $\frac{r+dk}{p+ck}$, see~\cref{f8}
  for an illustration.

  \begin{figure}[ht]
    \centering
    \includegraphics[width=1\textwidth]{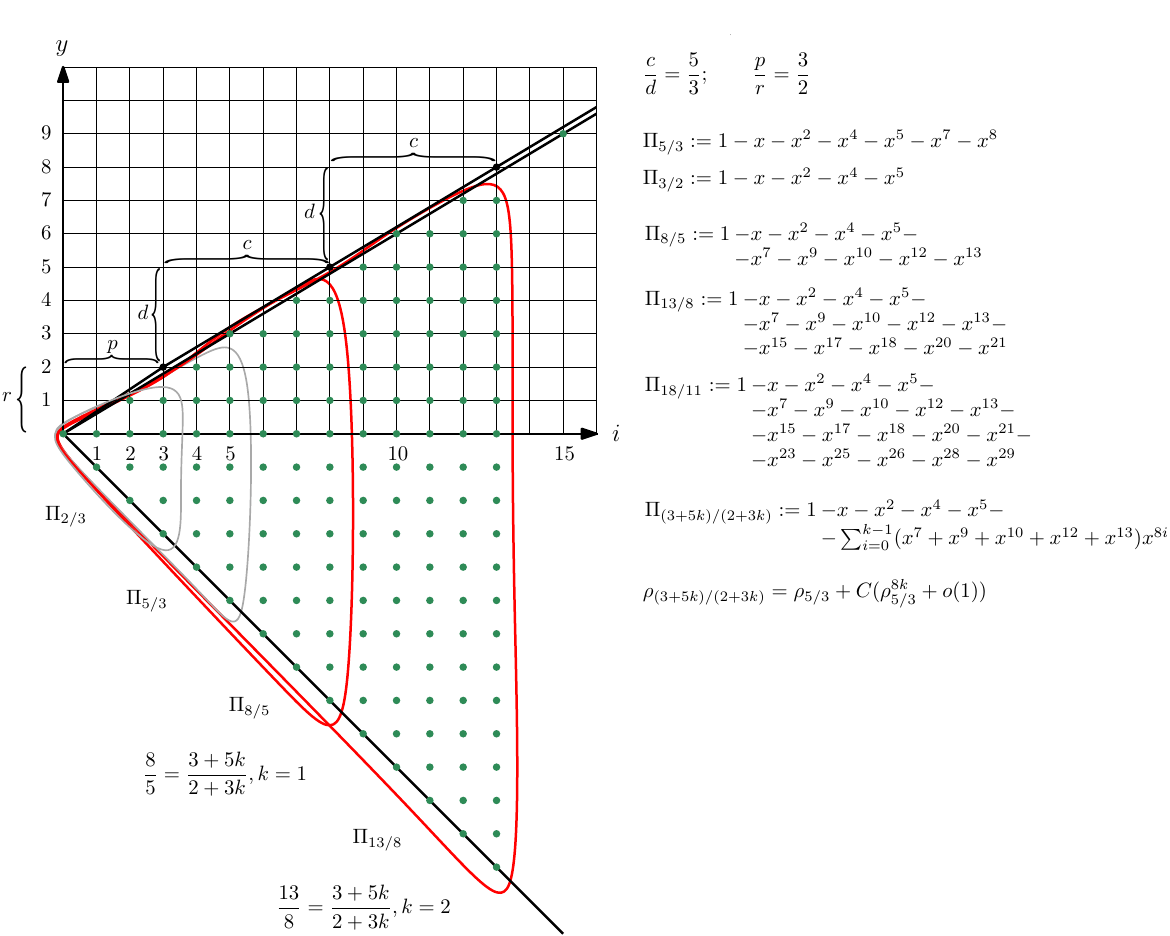}
    \caption{Geometric representation of polynomials $\Pi_{5/3},\Pi_{3/2}$
      and $\Pi_{(3+5k)/(2+3k)}$ for some values of $k\ge1$.}
    \label{f8}
  \end{figure}

  Note that from the construction of Stern--Brocot tree, $\frac{p}{r}
  < \frac{c}{d}$.  From the mediant inequality it follows that
  $\frac{p}{r} < \frac{p+ck}{r+dk} < \frac{c}{d}$, for any integer $k
  > 0$. Next, let us show that there are no integer points whose
  horizontal coordinate is equal to $i, 0 < i < p+ck$, lying strictly
  between the lines with slopes $\frac{r+dk}{p+ck}$ and $\frac{d}{c}$.
  Consider the triangle $ABC$ with following integer coordinates:
  $$
    {
    \everymath={\displaystyle}
    \begin{array}{l}
      A = (0,0);\\
      B = (c(k+1), d(k+1));\\
      C = (p+ck, r + dk).
    \end{array}
    }
    $$
    The area of $ABC$ is
    $$
    \frac{1}{2}
    \begin{vmatrix}
       c (k+1) & p + ck\\
       d (k+1) & r + dk
    \end{vmatrix} =
     \frac{1}{2} (k+1) \begin{vmatrix}
      c & p \\
      d & r
     \end{vmatrix}
     = \frac{k+1}{2} (cr - dp) = \frac{k+1}{2},
     $$
  because $\frac{p}{r} < \frac{c}{d}$ and $(cr - dp) = 1$ by a joli
  property of the Stern--Brocot tree (see~\cite{co}).  Pick's
  Theorem~\cite{pick} implies that the area of triangle $ABC$ is equal
  to $I + B/2 - 1$, where $I$ is the number of its interior points and
  exactly $B$ points lie on the boundary.  We have $B = k+3$, because
  $\frac{p+ck}{r+dk}$ and $\frac{c}{d}$ are irreducible fractions. We
  conclude that $I=0$ and there are no interior points inside
  $ABC$. This fact can be used to simplify $\Pi_{\frac{p+ck}{r+dk}}$
  from \cref{eq3} by considering only points lying on or below the
  line of slope $d/c$ and on or above the line of slope $(-c)$.
  Taking this into account, we rewrite the equation
  $\Pi_{\frac{p+ck}{r+dk}} = 0$, i.e.\ the equation
  $$
  1 - x^{p+r+(c+d)k} - \sum_{i=0}^{p+ck-1} x^{1 + i + \left\lfloor
        \frac{i (r+dk) }{p+ck} \right\rfloor} = 0
  $$
  by decomposing the internal sum representing the sum over all
  points in a triangle (see~\cref{f8}) by the sum over corresponding
  triangles and rectangles, which then allows us to simplify the sum
  by using the summation formula for geometric progressions:
  {
    \small
    \begin{align}
      & 1 - \sum_{j=0}^{p} x^{1+j+\left\lfloor   \frac{jd}{c} \right\rfloor}  -
       \sum_{j=p+1}^{p+c} x^{1+j+\left\lfloor   \frac{jd}{c} \right\rfloor}
       \sum_{i=0}^{k-1} x^{(c+d)i} = 0, \label{eq3.5}\\
       & 1 - \sum_{j=0}^{p} x^{1+j+\left\lfloor   \frac{jd}{c} \right\rfloor}  -
       \sum_{j=p+1}^{p+c} x^{1+j+\left\lfloor   \frac{jd}{c} \right\rfloor}
       \frac{1-x^{(c+d)k}}{1-x^{c+d}} = 0, \nonumber\\
       & 1 - \sum_{j=0}^{p} x^{1+j+\left\lfloor   \frac{jd}{c} \right\rfloor}
       - x^{c+d} + x^{c+d} \sum_{j=0}^{p} x^{1+j+\left\lfloor   \frac{jd}{c} \right\rfloor}
       - \sum_{j=p+1}^{p+c} x^{1+j+\left\lfloor   \frac{jd}{c} \right\rfloor}
       + x^{(c+d)k} \sum_{j=p+1}^{p+c} x^{1+j+\left\lfloor   \frac{jd}{c} \right\rfloor}
       = 0, \nonumber
    \end{align}
    \begin{align}
        & 1 - x^{c+d} - \sum_{j=0}^{p+c} x^{1+j+\left\lfloor   \frac{jd}{c} \right\rfloor}
         + \sum_{j=c}^{p+c} x^{1+j+d+\left\lfloor   \frac{(j-c)d}{c} \right\rfloor}
        + x^{(c+d)k} \sum_{j=p+1}^{p+c} x^{1+j+\left\lfloor   \frac{jd}{c} \right\rfloor}
        = 0, \nonumber\\
        & 1 - x^{c+d} - \sum_{j=0}^{c-1} x^{1+j+\left\lfloor   \frac{jd}{c} \right\rfloor}
        + x^{(c+d)k} \sum_{j=p+1}^{p+c} x^{1+j+\left\lfloor   \frac{jd}{c} \right\rfloor}
        = 0.  \label{eq4}
    \end{align}
  }%
  At some point, we multiply both sides by $(1-x^{c+d})$. All roots of
  Equation~\eqref{eq3.5} are also the roots of Equation~\eqref{eq4},
  but the latter is additionally satisfied by the roots of unity $1 =
  x^{c+d}$, the modulus of which is greater than
  $\rho_{\frac{p+ck}{r+dk}}$.

  Finally, since the first three terms of the sum are equal to
  $\Pi_{c/d}(x)$, we conclude that for $x = \rho_{\frac{p+ck}{r+dk}}$
  we have
  \begin{align}
        & - \Pi_{c/d}(x) = x^{(c+d)k} \sum_{j=p+1}^{p+c} x^{1+j+\left\lfloor   \frac{jd}{c} \right\rfloor}.\nonumber
    \end{align}
  Using the same method as in the proof of \cref{prop4}
  and denoting by $\Pi_{c/d}'$ the derivative of $\Pi_{c/d}$ we obtain
  $$
  \rho_{\frac{p+ck}{r+dk}} = \rho_{c/d} + C 
  \rho^{(c+d)k}_{c/d}\big( 1 + o(1) \big),
  $$
  where $$
  C = \frac{\sum_{j=p+1}^{p+c} \rho_{c/d}^{1+j+\left\lfloor   \frac{jd}{c} \right\rfloor}}{
    -\Pi_{c/d}'\left( \rho_{c/d} \right)}.
  $$
\end{proof}

\begin{proposition}
  For natural $k \ge 2$, the smallest in modulus root
  $\rho^+_{\frac{p+c(k-1)}{r+d(k-1)}}$ of
  $\Pi^+_{\frac{p+c(k-1)}{r+d(k-1)}}$ is
  $$\rho_{\frac{p+ck}{r+dk}} = \rho_{c/d} + C
  \rho^{(c+d)k}_{c/d} \big(1 + o(1) \big),$$ where $k\ge 1$,
  $\frac{p+c}{r+d}$ is the left child of $c/d$ in the Stern--Brocot tree,
  $C$ is the same constant as in \cref{prop6} not depending on $k$, and $\rho_{c/d}$ is the
  smallest in modulus root of $\Pi_{c/d}$.
  \label{prop7}
\end{proposition}

\begin{proof}
  Recall that \cref{proposition:jump} defines $\Pi^+_{a/b} := \Pi_{a/b} -
  x^{a+b} + x^{a+b+1}$. Having this in mind and adapting the equations
  from the proof of \cref{prop6} by writing $k-1$ in place of
  $k$ we obtain:
  $$
  {
    \everymath={\displaystyle}
    \begin{array}{l}
        1 - x^{c+d} - \sum_{j=0}^{c-1} x^{1+j+\left\lfloor   \frac{jd}{c} \right\rfloor}
        + x^{(c+d)(k-1)} \sum_{j=p+1}^{p+c} x^{1+j+\left\lfloor   \frac{jd}{c} \right\rfloor} - \\[1.4em]
        \hspace{2em} - x^{p+r+(c+d)(k-1)} (1-x^{c+d})
        + x^{p+r+(c+d)(k-1)+1} (1-x^{c+d}) = 0, \\[1.4em]
    \end{array}
  }
  $$
  \begin{equation}
  {
    \everymath={\displaystyle}
    \begin{array}{l}

        1 - x^{c+d} - \sum_{j=0}^{c-1} x^{1+j+\left\lfloor   \frac{jd}{c} \right\rfloor}
        + x^{(c+d)(k-1)} \times \\[1.4em]
        \hspace{2em} \times \left(
        \sum_{j=p+1}^{p+c} x^{1+j+\left\lfloor   \frac{jd}{c} \right\rfloor}
        - x^{p+r} + x^{p+r+(c+d)} + x^{p+r+1} - x^{p+r+(c+d)+1}
        \right) = 0.\\[1.4em]
        \label{eq8}
    \end{array}
  }
  \end{equation}
  Now, using the fact that $cr - dp = 1$ and a kind of geometrical
  argument (as in the previous proof), it can be shown that
  $${
    \everymath={\displaystyle}
    \begin{array}{l}
  \sum_{j=p+1}^{p+c} x^{1+j+\left\lfloor   \frac{jd}{c} \right\rfloor}
  - x^{p+r} + x^{p+r+(c+d)} + x^{p+r+1} - x^{p+r+(c+d)+1}
  - x^{(c+d)} \sum_{j=p+1}^{p+c} x^{1+j+\left\lfloor   \frac{jd}{c} \right\rfloor} = \\[1.4em]
  \hspace{2em} =  \left( 1-x^{c+d} - \sum_{j=0}^{c-1} x^{1+j + \left\lfloor \frac{jd}{c} \right\rfloor} \right) \left( x^{c+d+p+r} - x^{p+r} \right).
    \end{array}
  }
  $$ Adapting the techniques from the proofs of \cref{prop4,prop5}
  we see that the smallest in modulus root
  of~\eqref{eq8} is equal to the smallest in
  modulus root of
  $$
  {
    \everymath={\displaystyle}
    \begin{array}{l}
      1 - x^{c+d} - \sum_{j=0}^{c-1} x^{1+j+\left\lfloor   \frac{jd}{c} \right\rfloor}
      + x^{(c+d)(k-1)} \left(
       x^{(c+d)} \sum_{j=p+1}^{p+c} x^{1+j+\left\lfloor   \frac{jd}{c} \right\rfloor}
       + o(1)
      \right) = 0,\\[1.4em]
    \end{array}
  }
  $$
  where $o(1)$ is considered as $k \to \infty$.
  Comparing it to \cref{eq4} from the proof of \cref{prop6},
  we see that the claimed result follows.
\end{proof}

\subsection{Summary of fractal results and an open question}

To summarize the results of the previous two propositions, let $q$
denote $\frac{c}{d}$ and let $q^\ast_k$ denote $\frac{p+ck}{r+dk}$,
i.e.\ the $k$th left approximation of $q$ in the Stern--Brocot tree.
From \cref{prop6,prop7} the following result follows.

\begin{proposition}
The graph of the function \( \Phi(q) \) rescaled on the intervals
$(q_k^\ast, q]$ tends to a constant function on the semi-open
  intervals \( (q_{k-1}^\ast, q_k^\ast] \):
\[
  \lim_{k \to \infty}
  \dfrac
    {\Phi(q) - \lim_{\delta \to 0^+} \Phi(q_{k-1}^\ast + \delta)}
    {\Phi(q) - \Phi(q_k^\ast)} = 1.
\]
The ratio between consecutive constants on the rescaled graph tends to \( \rho_q^{c+d} \):
\[
  \lim_{k \to \infty}
  \dfrac
    {\Phi(q) - \Phi(q_k)}
    {\Phi(q) - \Phi(q_{k-1})} = \rho_{q}^{c+d}.
    \]
\end{proposition}

Using a similar technique, it is possible to demonstrate a resembling
picture using the right child of $1$ in the Stern--Brocot tree,
generating intervals of the form \( [1, \frac{k+1}{k}) \), and, more
  generally, using the right approximation of any rational point \(
  q=c/d \).

\medskip
  
From \cref{proposition:jump} it follows that the function $\Phi(q)$ has the
highest jump at point $q=1$, the second highest jump is at $q=1/2$, the
third highest jump appears when $q=2$, see \cref{f2}.  The sequence of positive
rational numbers ordered by corresponding jumps of the function
$\Phi(q) = \lim_{n \to \infty} |\W_{q,n+1}| / |\W_{q,n}|$ starts with
  $$
  1, \frac{1}{2}, 2, \frac{1}{3}, \frac{1}{4}, 3, \frac{2}{3}, \frac{1}{5}, \frac{1}{6}, \frac{3}{2}, \frac{1}{7}, 4,
  \frac{2}{5}, \frac{1}{8}, \frac{1}{9}, \frac{1}{10}, \frac{3}{4}, \frac{1}{11}, \frac{2}{7}, \frac{1}{12}, 5,
  \frac{3}{5}, \frac{1}{13}, \frac{4}{3}, \frac{1}{14}, \frac{2}{9},
  \frac{1}{15}, \frac{1}{16}, \frac{5}{2}, \frac{1}{17},...
  % \frac{3}{7}, \frac{1}{18},
%  \frac{2}{11}, \frac{1}{19}, 6, \frac{4}{5}, \frac{1}{20}, \frac{1}{21}, \frac{3}{8}, \frac{5}{3}, \frac{1}{22}, \frac{2}{13}, \frac{1}{23}, \frac{1}{24}, \frac{1}{25}, \frac{1}{26}, \frac{2}{15}, \frac{3}{10},
%   \frac{1}{27}, \frac{1}{28}, \frac{5}{4}, \frac{1}{29}, 7, ..
  $$

  \begin{question}
    Is it possible to explain this sequence without polynomial root calculations?
  \end{question}

\subsection*{Acknowledgements}

Sergey Dovgal was supported by the EIPHI Graduate School (contract
ANR-17-EURE-0002), FEDER R\'{e}gion Bourgogne Franche-Comt\'{e}, and Sergey
Kirgizov was supported in part by the project ANER ARTICO funded by
Bourgogne-Franche-Comt\'{e} region (France) and ANR-22-CE48-0002 funded by
l'Agence Nationale de la Recherche.  Special thanks to Michael A.
Allen, Jean-Luc Baril, Khaydar Nurligareev, Nathana\"{e}l Hassler and
Vincent Vajnovszki for inspiring and helpful discussions, to Vlady
Ravelomanana for helping to connect the two authors, and to Natalia
Kharchenko for helping to improve the current manuscript.  We thank
the anonymous reviewers who provided very good comments.

%BIBLIOGRAPHY
% You do not have to use the same format for your references, but 
%    include everything in this file.
% If you use BibTeX to create a bibliography, copy the .bbl file into here.
% We recommend you use \doi{...} and \arxiv{...} like the examples below,
% as they give a short display form with an active link to the full url.


\begin{thebibliography}{99}

\bibitem{auto}
  Allouche, J.-P., and Shallit, J.
  {\em Automatic Sequences: Theory, Applications, Generalizations}.
  Cambridge University Press, 2003

\bibitem{cointossing}
  Bacher, A., Bodini, O., Hwang, H. K., and Tsai, T. H.
  {\em Generating Random Permutations by Coin Tossing: Classical
        Algorithms, New Analysis, and Modern Implementation}.
  ACM Transactions on Algorithms, 13(2), 2017, pp. 1--43.

\bibitem{barcucci_D}
  Barcucci, E., Bernini, A., Bilotta, S., Pinzani R.
  {\em Dyck Paths Enumerated by the $\Q$-bonacci Numbers}.
  GASCom 2024, in Electronic Proceedings in Theoretical Computer Science, 403,
  2024, pp. 49-53.
  
\bibitem{barcucci_R}
  Barcucci, E., Bernini, A., Bilotta, S., Pinzani R.
  {\em Rational Dyck paths}.
  Journal of Integer Sequences, 28(3), 2025, Article 25.3.2.
  
\bibitem{ourfibo}
  Baril, J.-L., Kirgizov, S. and Vajnovszki, V.
  {\em Gray codes for Fibonacci q-decreasing words}.
  Theoretical Computer Science, 927, 2022, pp. 120--132.

\bibitem{be} Bernoulli, J. III.  {\em Sur une nouvelle espece de
  calcul}.  Recueil pour les Astronomes, 1, Berlin, 1771, pp. 255--284.

\bibitem{ber} Berstel, J., Lauve, A., Reutenauer, C. and Franco, S.
  {\em Combinatorics on words: Christoffel words and repetitions in words}.
  CRM Monograph Series, 27, American Mathematical Society, 2009.

\bibitem{brocot}
  Brocot, A.
  {\em Calcul des rouages par approximation, nouvelle m\'{e}thode}.
  Revue Chronom\'{e}trique, 3, 1861, pp. 186--194.

\bibitem{cri} Christoffel, E. B.
  {\em Observatio arithmetica}.
  Annali di Matematica Pura ed Applicata, 6, 1875, pp. 145--152.

\bibitem{denjoy} Denjoy, A.
  {\em Sur une fonction r{\'e}elle de {Minkowski}}.
  Journal de Math{\'e}matiques Pures et Appliqu{\'e}es. Neuvi{\`e}me S{\'e}rie,
  17, 1938, pp 105--151.

\bibitem{ei}  
  Eğecioğlu, \"O., and Iršič, V.
  {\em Fibonacci-run graphs {I}: Basic properties}.
  Discrete Applied Mathematics, 295, 2021, pp. 70--84.
  
\bibitem{fei}
  Feinberg, M.
  {\em Fibonacci-Tribonacci}.
  Fibonacci Quarterly, 1(3), 1963, pp. 71--74.

\bibitem{fla}
  Flajolet, P., and Sedgewick, R.
  {\em Analytic Combinatorics}.
  Cambridge University Press, 2009.

\bibitem{co}
  Graham, R.L., Knuth, D.E. and Patashnik, O.
  {\em Concrete Mathematics: A Foundation for Computer Science} 2nd ed.,
  Addison-Wesley, 1994

\bibitem{hassler}
  Hassler, N., Vajnovszki, V. and Wong, D.
  {\em Greedy gray codes for some restricted classes of binary words}.
  In Brlek, S. and Ferrari, L. (eds) Proc. of the 13th ed.
  of the conf. on Random Generation of Combinatorial Structures (GASCom), Polyominoes and Tilings, Bordeaux, France, 24-28th June 2024,
  Electronic Proceedings in Theoretical Computer Science, 403, Open Publishing Association, 2024, pp. 108--112.
  
\bibitem{ki}
  Kirgizov, S.
  {\em $\Q$-bonacci words and numbers}.
  Fibonacci Quarterly, 60(5), 2022, pp. 187--195.

 \bibitem{knuth3}
   Knuth, D. E.
   {\em The Art of Computer Programming, Volume 3: Sorting and Searching} 2nd ed.,
   Addison-Wesley, 1998.

\bibitem{lo} Lothaire, M.
  {\em Algebraic Combinatorics on Words}.
  Cambridge University Press, 2002.

\bibitem{miles}
  Miles Jr., E. P.
  {\em Generalized Fibonacci numbers and associated matrices}.
  The American Mathematical Monthly, 67(8), 1960, pp. 745--752.

\bibitem{minkowski}
  Minkowski, H.
  {\em Zur Geometrie der Zahlen}.
  Verhandlungen des III. internationalen Mathematiker-Kongresses in Heidelberg, Berlin, 1904, pp. 164--173.

\bibitem{orlov}
  Orlov, A. G.
  {\em On asymptotic behavior of the Taylor coefficients of algebraic functions}.
  Siberian Mathematical Journal, 25(5), 1994, pp. 1002--1013.
  
\bibitem{pick}
  Pick, G.
  {\em Geometrisches zur Zahlenlehre}.
  Sitzungsberichte des deutschen naturwissenschaftlich-medicinischen
  Vereines f\"{u}r B\"{o}hmen "Lotos" in Prag. (Neue Folge), 19, 1899, pp. 311--319.

\bibitem{slides-sedgewick}
  Sedgewick, R.
  {\em Analytic Combinatorics. Online Course. Chapter 4: Complex Analysis, Rational and Meromorphic Asymptotics}.
  \url{https://ac.cs.princeton.edu/online/slides/AC04-Poles.pdf},
  Princeton University, 2024.

\bibitem{stern}
  Stern, M. A.
  {\em Ueber eine zahlentheoretische Funktion}.
  Journal f\"{u}r die reine und angewandte Mathematik, 55, 1858, pp. 193--220.

\bibitem{wang}
  Wong, D., Liu, B.,  Lam, C.-T., and Im, M.
  {\em Generating cyclic 2-Gray codes for Fibonacci q-decreasing words}.
  In Uehara, R., Yamanaka, K., Yen, H.-C. (eds)
  WALCOM: Algorithms and Computation. WALCOM 2024, Kanazawa, Japan, 18-20th March 2024,
  Lecture Notes in Computer Science, vol. 14549, Springer, Singapore, 2024, pp. 91--102.
\end{thebibliography}
\end{document}